\documentclass[12pt]{amsproc}
\usepackage[cp1251]{inputenc}
\usepackage[T2A]{fontenc} 
\usepackage[russian]{babel} 
\usepackage{amsbsy,amsmath,amssymb,amsthm,graphicx,color} 
\usepackage{amscd}
\def\le{\leqslant}
\def\ge{\geqslant}

\newtheorem{prop}{Предложение}

\theoremstyle{definition}

\theoremstyle{remark}

\begin {document}
\unitlength=1mm
\title[Обобщенная формула Гассерта-Шора]
{Обобщённая формула Гассерта-Шора как комбинаторная теорема Стокса}
\author{Г. Г. Ильюта}
\email{ilgena@rambler.ru}
\address{}
\thanks{Работа поддержана грантом РФФИ-20-01-00579}


\begin{abstract}
Мы обобщим формулу Гассерта-Шора для числовых полугрупп.

We generalize Gassert-Shor formula for numerical semigroups.
\end{abstract}

\maketitle

  1. Введение. Множество $S\subset \mathbb Z_{\ge 0}$ называется числовой полугруппой, если оно содержит $0$, замкнуто относительно сложения и имеет конечное дополнение в $\mathbb Z_{\ge 0}$ \cite{5}, \cite{6}. Далее предполагаем, что $S\neq \mathbb Z_{\ge 0}$. Пусть
$$
C(S):=\mathbb Z_{\ge 0}\setminus S=\{c_1,\dots,c_n|c_1<\dots<c_n\},
$$
$F(S):=\max C(S)=c_n$ -- число Фробениуса полугруппы $S$. Для всей статьи мы фиксируем ненулевое $m\in S$ и для $i\in I_m:=\{0,1,\dots,m-1\}$ обозначим через $A_i$ количество элементов множества $C(S)$, которые сравнимы с $i$ по $\mod m$. Множество Апери $S^{(m)}$ определяется следующим образом:
$$
S^{(m)}:=\{s\in S:s-m\notin S\}.
$$
Другими словами, $S^{(m)}$ состоит из минимальных элементов в пересечениях полугруппы $S$ с классами вычетов по $\mod m$ -- для каждого $i\in I_m$
$$
\{i,i+m,\dots,i+m(A_i-1)\}\subset C(S),            \eqno (1)
$$
если $A_i>0$. В частности, $S^{(m)}$ образует содержащую $0$ полную систему вычетов по $\mod m$. Эти факты записаны в \cite{5} в виде равенства
$$
S=S^{(m)}+m\mathbb Z_{\ge 0}.                      \eqno (2)
$$
Приведём пример множества Апери: $\mathbb Z_{\ge 0}^{(m)}=I_m$. Полагая
$$
S^{(m)}=\{0=a_0,a_1,\dots,a_{m-1}:a_i\equiv i\mod m,i\in I_m\},
$$
имеем равенство $a_i=mA_i+i$ для всех $i\in I_m$.

  Дополнение $C(S)$ в $\mathbb Z_{\ge 0}$ и множество Апери $S^{(m)}$ являются двумя способами описания числовой полугруппы $S$. В \cite{1} доказано равенство, связывающее эти два определения -- своего рода формула перехода от одного базиса к другому: для произвольной функции $f$ на $\mathbb Z_{\ge 0}$
$$
\sum_{i\in C(S)}(f(i+m)-f(i))=\sum_{i=0}^{m-1}(f(a_i)-f(i)). \eqno (3)
$$
Из равенства (2) следует более короткое, чем в \cite{1}, доказательство формулы (3): суммируем в пересечении множества $C(S)$ с классом вычетов (1)
$$
(f(i+m)-f(i))+(f(i+2m)-f(i+m))+\dots+(f(i+A_im)-f(i+(A_i-1)m))
$$
$$
=f(i+A_im)-f(i)=f(a_i)-f(i)
$$
и затем суммируем по классам вычетов. Доказательство более общих фактов в п. 3 отличается только порядком суммирования. 

  Мы получим некоторые обобщения формулы (3) (будем называть её формулой Гассерта-Шора), отвечающие путям в дереве числовых полугрупп, а точнее простым разностным соотношениям для произвольной функции $H(S)$ на вершинах пути.  Для числовой полугруппы $S$ множество $\bar S:=S\cup F(S)$ также является числовой полугруппой \cite{6}, p. 33. Дерево числовых полугрупп определяется в \cite{6}, p. 92, как граф, вершины которого отвечают числовым полугруппам, а рёбра -- парам $(S,\bar S)$. Простейшее разностное соотношение -- это представление разности значений функции $H$ на концах пути в виде суммы разностей значений функции $H$ на концах входящих в путь рёбер. Формула (3) является частным случаем этого соотношения (см. п. 2) для канонического (кратчайшего) пути в дереве числовых полугрупп, соединяющего отвечающую полугруппе $S$ вершину с корнем дерева, отвечающим полугруппе $\mathbb Z_{\ge 0}$. Вершины этого пути (кроме корня) взаимно однозначно отвечают элементам множества $C(S)$ -- числам Фробениуса отвечающих вершинам числовых полугрупп. В п. 2 приводятся более общие разностные соотношения, в которых наборам последовательных вершин пути сопоставлены некоторые определители, но для краткости и наглядности примеры в п. 3 и п. 4 мы приведём только для простейшего разностного соотношения и его мультипликативного варианта. Представление формулы (3) как разностного соотношения подсказывает аналогию с теоремой Стокса (или Ньютона-Лейбница): интеграл по границе области равен интегралу от дифференциала по самой области. Считаем множество $C(S)$ областью, по которой идёт дискретное интегрирование, множество Апери $S^{(m)}$ -- границей этой области, а разность значений функции $H$ на отвечающих вершинам ребра числовых полугруппах (в частности, разность $f(i+m)-f(i)$ в формуле (3)) -- дискретным дифференциалом.

  Связь дерева числовых полугрупп с множествами Апери основана на следующем замечании: множества Апери $S^{(m)}$ и $\bar S^{(m)}$ отличаются одним элементом (см. п. 2). В п. 3 мы рассмотрим примеры функций $H(S)$, которые являются функциями от $S^{(m)}$ и для которых соответствующие разностные соотношения наиболее близки к формуле (3), в частности, каждое слагаемое в левых частях этих соотношений имеет множитель вида $f(i+m)-f(i)$, где $i\in C(S)$. Эти множители интерпретированы как разности переменных в хорошо известных рекуррентных соотношениях для элементарных и полных симметрических функций (или их производящих функций), а также для разделённых разностей Ньютона. Эти рекуррентные соотношения имеют одно общее свойство -- в них присутствуют пары наборов переменных, отличающиеся одним элементом. Также в п. 3 рассматривается пример функции $f(i)$, связанный с обобщениями классических многочленов Бернулли. Разностные соотношения для этих обобщений \cite{3} хорошо согласованы с нашим дискретным дифференциалом.

   Наш подход к обобщению формулы Гассерта-Шора (3) аналогичен подходу И. Макдональда к формуле Соломона для многочлена Пуанкаре $G(q)$ конечной группы $G$, порождённой отражениями, \cite{4}. Эта аналогия прослеживается в п. 5. Там же приводятся ещё два равенства из \cite{2} и \cite{7}, которые, на наш взгляд, похожи на теорему Стокса и на формулу Гассерта-Шора (3).

\bigskip

  2. Дерево числовых полугрупп и множества Апери. При добавлении к полугруппе $S$ её числа Фробениуса $F(S)$ меняется пересечение полугруппы $S$ только с одним классом вычетов по $\mod m$, а значит в множестве $S^{(m)}$ меняется только один элемент. Из сказанного в первом абзаце Введения следует, что
$$
\bar S^{(m)}\setminus\{F(S)\}=S^{(m)}\setminus\{F(S)+m\},
$$
другими словами, множество Апери $\bar S^{(m)}$ получается из множества Апери $S^{(m)}$ заменой числа $F(S)+m$ на число $F(S)$. Действительно, $F(S)+m\in S^{(m)}$, так как $F(S):=\max C(S)\notin S$, а следующее рассуждение показывает, что $F(S)\in \bar S^{(m)}$: $F(S)-m\notin \bar S$, так как в противном случае $F(S)-m\in S$, а значит $m+(F(S)-m)=F(S)\in S$. 
 
  Имея для $l>0$ конечную последовательность числовых полугрупп, в которой соседние полугруппы отличаются одна от другой добавлением числа Фробениуса (путь в дереве числовых полугрупп)
$$
V_0-V_1-\dots-V_{l-1}-V_l,
$$
и функцию $H(V)$ на множестве $\{V_0,\dots,V_l\}$, напишем для $0<p\le l$ равенство для суммы определителей, отличающихся одной строкой,
$$
\sum_{i=0}^{l-p}D(H(V_{i+1})-H(V_i),\dots,H(V_{i+p})-H(V_{i+p-1}))$$
$$
=D(H(V_{l-p+1})-H(V_0),\dots,H(V_l)-H(V_{p-1})),
$$
где $D(y_1,\dots,y_p)$ -- определитель, первая строка которого совпадает с $y_1 \dots y_p$, а остальные элементы (комплексные числа или коммутирующие переменные) фиксированы. Этот определитель можно считать формальной линейной комбинацией элементов первой строки и поэтому в качестве функции $H(V)$ допускать функцию с произвольной областью значений. Для $p=1$ имеем равенство
$$
\sum_{i=0}^{l-1}(H(V_{i+1})-H(V_i))=H(V_l)-H(V_0)   \eqno (4)
$$
или же в мультипликативном варианте (если дроби определены)
$$
\prod_{i=0}^{l-1}\frac{H(V_{i+1})}{H(V_i)}
=\frac{H(V_l)}{H(V_0)}.                             
$$

  В оставшейся части статьи мы будем рассматривать для произвольной числовой полугруппы $S$ только её канонический путь в дереве числовых полугрупп
$$
S=S_n-S_{n-1}-\dots-S_1-S_0=\mathbb Z_{\ge 0},
$$
где $S_i:=\mathbb Z_{\ge 0}\setminus\{c_1,\dots,c_i\}$, $i=1,\dots,n$. Для каждого $i=1,\dots,n$ $\bar S_i=S_{i-1}$, $F(S_i)=c_i$ и 
$$
T_i:=S_i^{(m)}\setminus\{c_i+m\}=S_{i-1}^{(m)}\setminus\{c_i\},
$$
другими словами, множество Апери $S_{i-1}^{(m)}$ получается из множества Апери $S_i^{(m)}$ заменой числа $c_i+m$ на число $c_i$. Для канонического пути равенство (4) имеет вид
$$
\sum_{i=0}^{n-1}(H(S_{i+1})-H(S_i))
=H(S)-H(\mathbb Z_{\ge 0}).                    \eqno (5)
$$
Полагая в этом равенстве $H(S_i)=\sum_{a\in S_i^{(m)}}a$, получим формулу (3). В п. 4 приведён пример специализации для мультипликативного варианта этого разностного соотношения.

\bigskip

  3. Примеры функций на вершинах канонического пути. Введём обозначение для множества значений функции $f$ на элементах множества $A$
$$
f\{A\}:=\{f(a):a\in A\}.
$$
Для $0<p\le m$ и симметрической функции $F(x_1,\dots,x_p)$ пусть
$$
H(S_i):=\sum_{I\subset S_i^{(m)},|I|=p}F(f\{I\}).
$$
Для такой функции $H(S_i)$ формула (5) примет вид
$$
\sum_{i=1}^n\left(\sum_{J\subset T_i,|J|=p-1}F(f\{J\},f(c_i+m))-\sum_{J\subset T_i,|J|=p-1}F(f\{J\},f(c_i))\right)               
$$
$$
=\sum_{I\subset S^{(m)},|I|=p}F(f\{I\})-\sum_{I\subset I_m,|I|=p}F(f\{I\}),                               \eqno (6)
$$
в частности, для $p=m$ 
$$
\sum_{i=1}^n(F(f\{T_i\},f(c_i+m))-F(f\{T_i\},f(c_i)))
$$
$$
=F(f\{S^{(m)}\})-F(f\{I_m\}).
$$

  Полагая в формуле (6) $F(x_1,\dots,x_p)=(z+x_1)\dots(z+x_p)$, получим

\begin{prop}\label{prop1} 
$$
\sum_{i=1}^n(f(c_i+m)-f(c_i))\left(\sum_{J\subset T_i,|J|=p-1}\prod_{a\in J}(z+f(a))\right)
$$
$$
=\sum_{I\subset S^{(m)},|I|=p}\prod_{a\in I}(z+f(a))-\sum_{I\subset I_m,|I|=p}\prod_{a\in I}(z+f(a)),                                                       \eqno (7)
$$
в частности, для $p=m$
$$
\sum_{i=1}^n(f(c_i+m)-f(c_i))\prod_{a\in T_i}(z+f(a))
$$
$$
=\prod_{a\in S^{(m)}}(z+f(a))-\prod_{a\in I_m}(z+f(a)).
$$
\end{prop}

  Полагая в формуле (6) $F(x_1,\dots,x_p)=(z-x_1)^{-1}\dots(z-x_p)^{-1}$, получим

\begin{prop}\label{prop2} 
$$
\sum_{i=1}^n\frac{f(c_i+m)-f(c_i)}{(z-f(c_i))(z-f(c_i+m))}\left(\sum_{J\subset T_i,|J|=p-1}\prod_{a\in J}(z-f(a))^{-1}\right)
$$
$$
=\sum_{I\subset S^{(m)},|I|=p}\prod_{a\in I}(z-f(a))^{-1}-\sum_{I\subset I_m,|I|=p}\prod_{a\in I}(z-f(a))^{-1},                                   \eqno (8)
$$
в частности, для $p=m$
$$
\sum_{i=1}^n\frac{f(c_i+m)-f(c_i)}{(z-f(c_i))(z-f(c_i+m))}\prod_{a\in T_i}(z-f(a))^{-1}
$$
$$
=\prod_{a\in S^{(m)}}(z-f(a))^{-1}-\prod_{a\in I_m}(z-f(a))^{-1}.
$$
\end{prop}

  Элементарные и полные симметрические функции определяются равенствами
$$
e_k(x_1,\dots,x_p):=\sum_{i_1<\dots<i_k}x_{i_1}\dots x_{i_k},
$$
$$
h_k(x_1,\dots,x_p):=\sum_{i_1\le \dots\le i_k}x_{i_1}\dots x_{i_k},
$$
в частности,
$$
e_0(x_1,\dots,x_p)=h_0(x_1,\dots,x_p)=1,
$$
$$
e_1(x_1,\dots,x_p)=h_1(x_1,\dots,x_p)=\sum_{i=1}^px_i.
$$

\begin{prop}\label{prop3} Для $k=1,\dots,p$
$$
\sum_{i=1}^n(f(c_i+m)-f(c_i))\left(\sum_{J\subset T_i,|J|=p-1}e_{k-1}(f\{J\})\right)
$$
$$
=\sum_{I\subset S^{(m)},|I|=p}e_k(f\{I\})-\sum_{I\subset I_m,|I|=p}e_k(f\{I\}).
$$
В частности, для $p=m$ 
$$
\sum_{i=1}^n(f(c_i+m)-f(c_i))e_{k-1}(f\{T_i\})
=e_k(f\{S^{(m)}\})-e_k(f\{I_m\}),
$$
для $k=1$ получаем формулу (3).
\end{prop}

Доказательство. Полагаем в формуле (6) $F(x_1,\dots,x_p)=e_k(x_1,\dots,x_p)$ и используем соотношение
$$
e_k(x_2,\dots,x_p)-e_k(x_1,\dots,x_{p-1})=(x_p-x_1)e_k(x_2,\dots,x_{p-1}).
$$
Также можно использовать равенство
$$
(z+x_1)\dots(z+x_p)=\sum_{k=0}^pe_k(x_1,\dots,x_p)z^{p-k}
$$
в формуле (7) и сравнить коэффициенты при одинаковых степенях переменной $z$.

\begin{prop}\label{prop4} Для $k\ge 1$
$$
\sum_{i=1}^n(f(c_i+m)-f(c_i))\left(\sum_{J\subset T_i,|J|=p-1}h_{k-1}(f\{J\},f(c_i),f(c_i+m))\right)
$$
$$
=\sum_{I\subset S^{(m)},|I|=p}h_k(f\{I\})-\sum_{I\subset I_m,|I|=p}h_k(f\{I\}).
$$
В частности, для $p=m$ 
$$
\sum_{i=1}^n(f(c_i+m)-f(c_i))h_{k-1}(f\{T_i\},f(c_i),f(c_i+m))
=h_k(f\{S^{(m)}\})-h_k(f\{I_m\}),
$$
для $k=1$ получаем формулу (3).
\end{prop}

Доказательство. Полагаем в формуле (6) $F(x_1,\dots,x_p)=h_k(x_1,\dots,x_p)$ и используем соотношение
$$
h_k(x_2,\dots,x_p)-h_k(x_1,\dots,x_{p-1})=(x_p-x_1)h_{k-1}(x_1,\dots,x_p).
$$
Также можно использовать равенство
$$
z^p(z-x_1)^{-1}\dots(z-x_p)^{-1}=\sum_{k\ge 0}h_k(x_1,\dots,x_p)z^{-k}
$$
в формуле (8) и сравнить коэффициенты при одинаковых степенях переменной $z$.

  Определим разделённые разности Ньютона $\Delta(w(z);x_1,\dots,x_p)$ функции $w(z)$ как контурный интеграл (мы пропускаем общеизвестные детали)
$$
\Delta(w(z);x_1,\dots,x_p):=\frac{1}{2\pi i}\int\frac{w(z)dz}{(z-x_1)\dots(z-x_p)}.
$$
Умножая равенства в Предложении 2 на $w(z)/2\pi i$ и интегрируя, получим

\begin{prop}\label{prop5} 
$$
\sum_{i=1}^n(f(c_i+m)-f(c_i))\sum_{J\subset T_i,|J|=p-1}\Delta(w(z);f\{J\},f(c_i),f(c_i+m))
$$
$$
=\sum_{I\subset S^{(m)},|I|=p}\Delta(w(z);f\{I\})-\sum_{I\subset I_m,|I|=p}\Delta(w(z);f\{I\}),
$$
в частности, для $p=m$
$$
\sum_{i=1}^n(f(c_i+m)-f(c_i))\Delta(w(z);f\{T_i\},f(c_i),f(c_i+m))
$$
$$
=\Delta(w(z);f\{S^{(m)}\})-\Delta(w(z);f\{I_m\}),
$$
для $w(z)=z^m$ получаем формулу (3).
\end{prop}

  Приведём пример функции $f(i)$ в формуле (3) и в формулах из Предложений 1-5. Несколько обобщений многочленов Бернулли $B_n(t)$ объединены в следующем определении \cite{3}, р. 337: пусть для $[l]_q:=(q^l-1)/(q-1)$
$$
B_{n;q^l;y}^{(\alpha)}(t+1;\lambda):=\sum_{k=0}^n\binom{n}{k}q^{l(k-\alpha+1)y}[l]_q^kB_{k;q^l}^{(\alpha)}(t;\lambda),
$$
а $B_{k;q}^{(\alpha)}(t;\lambda)$ определяются производящей функцией 
$$
(-z)^\alpha \sum_{n=0}^\infty\frac{[\alpha]_q[\alpha+1]_q\dots [\alpha+n-1]_q}{[1]_q[2]_q\dots [n]_q}\lambda^nq^{n+t}e^{[n+t]_qz}=\sum_{n=0}^\infty B_{k;q}^{(\alpha)}(t;\lambda)\frac{z^n}{n!}.
$$
Для $n\ge 1$ имеем разностное соотношение \cite{3}, р. 338,
$$
\lambda q^{l(\alpha-1)}B_{n;q^l;y}^{(\alpha)}(t+1;\lambda)-B_{n;q^l;y}^{(\alpha)}(t;\lambda)=n[l]_qB_{n-1;q^l;y}^{(\alpha-1)}(t;\lambda).                                             \eqno (9)  
$$
Обозначим через $[x]$ целую часть $x\in\mathbb R$. Полагая
$$
f(i)=(\lambda q^{l(\alpha-1)})^{\left[\frac{x+i}{m}\right]}B_{n;q^l;y}^{(\alpha)}\left(\frac{i}{m};\lambda\right),
$$
и используя разностное соотношение (9), получим равенство
$$
f(i+m)-f(i)=(\lambda q^{l(\alpha-1)})^{\left[\frac{x+i}{m}\right]}n[l]_qB_{n-1;q^l;y}^{(\alpha-1)}\left(\frac{i}{m};\lambda\right),
$$
Поэтому формула (3) принимает вид
$$
n[l]_q\sum_{i\in C(S)}(\lambda q^{l(\alpha-1)})^{\left[\frac{x+i}{m}\right]}B_{n-1;q^l;y}^{(\alpha-1)}\left(\frac{i}{m};\lambda\right)
$$
$$
=\sum_{i=0}^{m-1}\left((\lambda q^{l(\alpha-1)})^{\left[\frac{x+a_i}{m}\right]}B_{n;q^l;y}^{(\alpha)}\left(\frac{a_i}{m};\lambda\right)-(\lambda q^{l(\alpha-1)})^{\left[\frac{x+i}{m}\right]}B_{n;q^l;y}^{(\alpha)}\left(\frac{i}{m};\lambda\right)\right).
$$
Аналогично преобразуются формулы из Предложений 1-5.

\bigskip

  4. Пример специализации для мультипликативного варианта разностного соотношения (5).

\begin{prop}\label{prop6} Для $k\ge 1$, $1\le p<m$
$$
\prod_{i=1}^n\prod_{J\subset T_i\setminus 0,|J|=p-1}\left(1+\frac{\sum_{j=1}^k\binom{k}{j}c_i^{k-j}m^j}{\sum_{a\in J}a^k+c_i^k}\right)=\frac{\prod_{I\subset S^{(m)}\setminus 0,|I|=p}\sum_{a\in I}a^k}{\prod_{I\subset I_m\setminus 0,|I|=p}\sum_{a\in I}a^k}.
$$
в частности, для $p=m-1$ 
$$
\prod_{i=1}^n\left(1+\frac{\sum_{j=1}^k\binom{k}{j}c_i^{k-j}m^j}{\sum_{a\in T_i}a^k+c_i^k}\right)
=\frac{\sum_{a\in S^{(m)}}a^k}{\sum_{a\in I_m}a^k},
$$
для $p=1$
$$
\prod_{i=1}^n\left(1+\frac{m}{c_i}\right)
=\frac{\prod_{a\in S^{(m)}\setminus 0}a}{(m-1)!}.
$$
\end{prop}

Доказательство. В мультипликативном варианте разностного соотношения (5) 
$$
\prod_{i=1}^n\frac{H(S_i)}{H(S_{i-1})}=\frac{H(S)}{H(\mathbb Z_{\ge 0})}.
$$
полагаем
$$
H(S_i):=\prod_{I\subset S_i^{(m)}\setminus 0,|I|=p}\sum_{a\in I}a^k
$$
Тогда
$$ 
\frac{H(S_i)}{H(S_{i-1})}=\prod_{J\subset T_i\setminus 0,|J|=p-1}\frac{\sum_{a\in J}a^k+(c_i+m)^k}{\sum_{a\in J}a^k+c_i^k}
$$
$$
=\prod_{J\subset T_i\setminus 0,|J|=p-1}\frac{\sum_{a\in J}a^k+c_i^k+\sum_{j=1}^k\binom{k}{j}c_i^{k-j}m^j}{\sum_{a\in J}a^k+c_i^k}
$$
$$
=\prod_{J\subset T_i\setminus 0,|J|=p-1}\left(1+\frac{\sum_{j=1}^k\binom{k}{j}c_i^{k-j}m^j}{\sum_{a\in J}a^k+c_i^k}\right).
$$

\bigskip

  5. Аналогия между формулой Макдональда и формулой Гассерта-Шора. Формулы Соломона и Макдональда для многочлена Пуанкаре $G(q)$ имеют вид
$$
G(q)=\prod_i\frac{q^{e_i+1}-1}{q-1}
$$
$$
=\prod_{r\in R_+(G)}\frac{q^{ht(r)+1}-1}{q^{ht(r)}-1},
$$
где $e_i$ показатели Кокстера группы $G$, $ht(r)$ -- высота (сумма координат в базисе из простых корней) корня $r$ из системы положительных корней $R_+(G)$ группы $G$ (одной из систем корней $A_n$, $B_n$, $C_n$, $D_n$, $E_6$, $E_7$, $E_8$, $F_4$, $G_2$). И. Макдональд получил второе равенство из первого с помощью мультипликативного варианта разностного соотношения
$$
\frac{q^{e_i+1}-1}{q-1}=\prod_{j=1}^{e_i}\frac{q^{j+1}-1}{q^j-1}
$$
и следующего известного факта: разбиение, составленное из показателей Кокстера $e_i$ сопряжено разбиению, составленному из количеств $b_i$ положительных корней высоты $i$,
$$
b_i:=|\{r\in R_+(G):ht(r)=i\}|.
$$
Этот факт можно представить как аналог формулы Гассерта-Шора (3): для произвольной функции $f$ на $\mathbb Z_{\ge 0}$
$$
\sum_{r\in R_+(G)}(f(ht(r)+1)-f(ht(r))]
=\sum_i(f(e_i+1)-f(1))                        \eqno (10)
$$
или же в мультипликативном варианте (если дроби определены)
$$
\prod_{r\in R_+(G)} \frac{f(ht(r)+1)}{f(ht(r))}=\prod_i\frac{f(e_i+1)}{f(1)}.
$$
Можно расширить аналогию и в обратном направлении, а именно, определим $m$-высоту вещественного числа $x$ как $ht_m(x):=[x/m]$ и пусть
$$
b_i(S):=|\{k\in C(S):ht_m(k)=i\}|.
$$
Из определений следует, что разбиение, составленное из чисел $A_i$, сопряжено разбиению, составленному из чисел $b_i(S)$. 

  Формулу (10) также можно рассматривать как комбинаторную теорему Стокса: областью дискретного интегрирования (аналогом дополнения $C(S)$ к числовой полугруппе $S$) считаем мультимножество высот корней из системы положительных корней $R_+(G)$ группы $G$ (или же само множество $R_+(G)$), границей этой области (аналогом множества Апери) -- набор чисел $d_i=e_i+1$, совпадающий с набором степеней базисных инвариантов группы $G$, а дискретным дифференциалом -- функцию
$$
f(ht(r)+1)-f(ht(r))=(f(ht(r)+1)-f(1))-(f(ht(r))-f(1)).
$$

  Полагая в формуле (10)
$$
f(i)=(\lambda q^{l(\alpha-1)})^iB_{n;q^l;y}^{(\alpha)}(i;\lambda),
$$
и используя разностное соотношение (9), получим для $n\ge 1$ следующее равенство
$$
n[l]_q\sum_{r\in R_+(G)}(\lambda q^{l(\alpha-1)})^{ht(r)}B_{n-1;q^l;y}^{(\alpha-1)}(ht(r);\lambda)
$$
$$
=\sum_i((\lambda q^{l(\alpha-1)})^{d_i}B_{n;q^l;y}^{(\alpha)}(d_i;\lambda)-\lambda q^{l(\alpha-1)}B_{n;q^l;y}^{(\alpha)}(1;\lambda)).
$$

  Приведём ещё два равенства, которые, на наш взгляд, похожи на теорему Стокса и на формулу Гассерта-Шора (3). Для $k\in \mathbb Z_{>0}$ пусть
$$
M(k):=\left\{i\in \mathbb Z_{>0}:\frac{k}{i}-\left[\frac{k}{i}\right]\ge\frac{1}{2}\right\}. 
$$
В \cite{2} появляется равенство
$$
\sum_{i=1}^{2k}\left(f(i)\left[\frac{2k}{i}\right]-2f(i)\left[\frac{k}{i}\right]\right)=\sum_{i\in M(k)}f(i), 
$$
вытекающее из соотношения
$$
[2x]-2[x]=\begin{cases}
          0,&\text{если $x-[x]<1/2$;}\\ 
          1,&\text{если $x-[x]\ge 1/2$.}
          \end{cases}
$$
В \cite{7} появляется равенство
$$
\sum_{i\ge 1}\left(f(i)\left[\frac{k}{i}\right]-f(i)\left[\frac{k-1}{i}\right]\right)=\sum_{i|k}f(i), 
$$
вытекающее из соотношения
$$
\left[\frac{k}{i}\right]-\left[\frac{k-1}{i}\right]
         =\begin{cases}
          0,&\text{если $i\nmid k$;}\\ 
          1,&\text{если $i|k$.}
          \end{cases}
$$

\end {document}